\documentclass[12pt]{amsart}
\usepackage{amssymb}
\usepackage{bm}

\usepackage[
width=31pc,
height=48pc,
margin=1in,
footskip=30pt,
]{geometry}
\usepackage{layout}

\usepackage{graphicx}
\usepackage{epsfig}
\usepackage{subcaption}
\usepackage{amsmath,amsfonts,epstopdf}
\usepackage{tikz-cd}
\usepackage{multirow}

\theoremstyle{plain}
\newtheorem{theorem}{Theorem}[section]
\newtheorem{lemma}[theorem]{Lemma}

\theoremstyle{definition}
\newtheorem{definition}[theorem]{Definition}

\numberwithin{equation}{section}

\renewcommand{\phi}{\varphi}

\title[A survey on the topological entropy of cubic polynomials]{A survey on the topological entropy of cubic polynomials}
\author{Noah Cockram and Ana Rodrigues}
\address{Department of Mathematics\\University of Exeter\\Exeter EX4 4QF, UK}
\begin{document}

\maketitle

\begin{abstract} In this paper we discuss two different existing algorithms for computing topological entropy and we perform one of them in order to compute the isentropes for cubic polynomials.

\end{abstract}

\section{Introduction}

Some results concerning topological entropy are known for maps of an interval. For instance, for  piecewise monotone maps of an interval, the topological entropy is given by
$$h(f)= \lim_{n \rightarrow \infty} \frac{1}{n} \log c_n,$$
where $c_n$ is the smallest number of intervals on which $f^n$ is monotone \cite{MS80}. We  say that $f$ has a constant slope $s$ if on each of its pieces of monotonicity the map is affine with  slope $s$. It is  known that for maps of this type  $h(f ) = max(0, \log s)$ \cite{MS80}.  Also if $f : I \rightarrow I$ is an interval map, then (see \cite{Mis79} and \cite{Mis80})
$\limsup_{n \rightarrow \infty} \frac{1}{n} \log Card \{x \in I : f^n(x) = x \} \geq h(f).$
A point $p$ is said to be {\it nonwandering} if for every neighbourhood of $p$ there is an integer $n>0$ such that $f^n(U) \cap U \neq \emptyset$ and the set of all nonwandering points for $f$ is called the {\it nonwandering set}. If the {\it nonwandering set} of $f$ is a finite number of periodic points, then $h(f)=0$.

When working on topological entropy of one-dimensional Dynamical Systems, one would like to have some simple method of computing it and that is in general a very difficult task. For instance, for piecewise strictly monotone interval maps, various computational methods are known  \cite{Blocketal}  but none of them is really general and simple. Even for Markov maps, before starting computations, which are relatively simple, because one computes only the spectral radius of a non-negative matrix \cite{Alsedaetal} one has to identify the Markov structure and find the transition matrix. This structure may be complicated and the matrix can be large.

Another important question when studying topological entropy is whether the level sets of topological entropy (sets for which the topological entropy is constant)  are locally connected. Very recently, Rempe and van Strien \cite{GS15}  proved that each isentrope within the space of trigonometric polynomials is connected.  Milnor and Tresser \cite{MT00} considered real cubic maps of the interval onto itself, either with positive or with negative leading coefficient proving that each locus of constant topological entropy in parameter space is a connected set.

Bruin and van Strien \cite{bruin2013} proved the following result regarding the structure of isentropes (level sets of constant entropy) of polynomial maps.

\begin{theorem}
\label{NonlocalConnect}
For $m \in \mathbb{N} \setminus \{0\}$, $\epsilon \in \{-1,1\}$, let $P_\epsilon^m$ be the space of $m$-modal polynomial maps on the interval $[-1,1]$, such that $f(-1) = \epsilon$ and $f(1) \in \{-1,1\}$.
For any $m \geq 3$, there exists a dense set $H \subset [0,\log(m - 1)]$ such that for each $h^* \in H$, the $h^*$-isentrope of $P_\epsilon^m$ is not locally connected.
\end{theorem}

The question of whether this theorem holds true for $m = 2$, i.e. for the space of \textit{boundary-anchored cubic polynomials}, is still an open problem.  

Many authors have proposed numerical algorithms for computing efficiently the topological entropy. For instance, in \cite{GB91}, an algorithm for piecewise monotone maps of an interval which are not necessarily continuous is presented, however, this algorithm proved not to be accurate. In \cite{BS97} an algorithm equivalent to the standard power method for finding eigenvalues of matrices (with a shift of the origin) in the cases for which the function is Markov is introduced.

We begin this survey by comparing the two more recent algorithms for computing topological entropy by Radulescu \cite{Anca2008} and one by  Amigo and Gimenez \cite{amigo2014}.
We then implement the algorithm in \cite{amigo2014} in order to study the topological entropy of boundary anchored cubic polynomials.

\section{Algorithms for computing the topological entropy}

In this section we 
present two different algorithms for computing topological entropy, one by Radulescu \cite{Anca2008} and one by  Amigo and Gimenez \cite{amigo2014}. 
We then compare the two algorithms.

\subsection{Radulescu's algorithm}

Now we will turn to a particular subset of quartic polynomials, in order to discuss the first of the two algorithmic approaches for computing entropy that we will investigate in this paper, Radulescu's algorithm \cite{Anca2008}.  

We will start by computing the topological entropy of quartic polynomial interval maps $f:[0,1] \rightarrow [0,1]$, such that $f$ is the composition of two logisitic maps, i.e.
$$
f_{\mu,\lambda} = q_\mu \circ q_\lambda,
$$
where, 
$$
q_a : [0,1] \rightarrow [0,1]; \hspace{3mm} x \mapsto ax(1 - x),
$$
for $a \in [0,4]$.  The space of such polynomials will be called $P^Q$, and is clearly parameterisable by two parameters $(\mu,\lambda) \in [0,4]^2$.  Thus, we can draw the isentropes of $P^Q$ in the square $[0,4]^2$, the coordinate $(x,y) \in [0,4]^2$ represents the polynomial with the parameters $(x,y)$.  Another space of interval maps we will be using is $P^T$, the space of sawtooth maps $S: [0,1] \rightarrow [0,1]$ of the form
$$
S = T_b \circ T_a,
$$
where $T_a$ is a tent map of the form
$$
T_a : [0,1] \rightarrow [0,1]; \hspace{3mm} x \mapsto a \bigg( \frac{1}{2} - \bigg| x - \frac{1}{2} \bigg| \bigg),
$$
for $a \in [0,2]$.
\par
The sawtooth map $T_a \circ T_b$ consists of line segments of slope of $\pm ab$, and thus has topological entropy equal to $\max\{0,\log ab = \log a + \log b\}$.  Thus, for our purposes, we will be considering sawtooth maps from the parameter space
\begin{align*}
\mathcal{S} &= \{(a,b) \in [0,2]^2 : \log a + \log b \geq 0\} \\
&= \{(a,b) \in [\tfrac{1}{2}, 2]^2, \log a + \log b \geq 0\}.
\end{align*}
Note that if $a < \frac{1}{2}$, i.e. $\log a < -\log 2$, then we must have $\log b > \log 2$, implying $b > 2$, which is a contradiction as $b \in [0,2]$.  Similarly $b < \frac{1}{2}$ implies $a > 2$, which is again a contradiction.  Hence, the two sets above are indeed equal.
\par
In order to compute the topological entropy of these maps, we will use the algorithm  \cite{Anca2008}.  This algorithm attempts to compute topological entropy by approximating it with sawtooth maps.  This algorithm exploits the defining property of these quartic maps, that they are a composition of two logistic maps.  It does this by considering the logistic maps in question individually, specifically looking at the behaviour of the critical point by alternating application of either map.  From this behaviour, and using a few results which Radulecu's paper proves, we have a theoretical foundation which uses the behaviour of these simpler logistic maps to determine the entropy of the larger quartic map.  

Let us now introduce {\it d'itineraries} (following \cite{Anca2008}).

Let $f_1 : I_1 \rightarrow I_2$, $f_2: I_2 \rightarrow I_1$ be quadratic maps.  We say that, for $x \in I_1$, the following sequence:
$$
x, f_1(x), f_2(f_1(x)), f_1(f_2(f_1(x))), ...
$$
is known as the \textit{diorbit} of $x$ under $(f_1,f_2)$.  If $x \in I_2$, then the diorbit of $x$ under $(f_1,f_2)$ is:
$$
x, f_2(x), f_1(f_2(x)), f_2(f_1(f_2(x))), ...
$$
Let $c_1 \in I_1$ be the critical point of $f_1$, and let $c_2 \in I_2$ be the critical point of $f_2$.  We define the alphabets $\{L_1,C_1,R_1\}$ and $\{L_2,C_2,R_2\}$ such that $L_i = [0,c_i)$, $C_i = \{c_i\}$, and $R_i = (c_i,1]$ for $i = 1,2$.  Let $(x_j)_{j\in \mathbb{N}}$ be the diorbit of $x$ under $(f_1,f_2)$.  The \textit{d'itinerary} of $x$ under $(f_1,f_2)$ is defined as the sequence $i(f_1,f_2)(x) = (A_j)_{j \in \mathbb{N}}$, where
\begin{equation*}
A_j = \Bigg\{ \begin{aligned}
                      &L_i \text{ if } x_j \in L_i \\
                      &C_i \text{ if } x_j = c_i \\
                      &R_i \text{ if } x_j \in R_i
                      \end{aligned}.
\end{equation*}
In a similar way to itineraries, a partial order can be defined on d'itineraries, even on d'itineraries under different maps!  This partial order is undefined in some cases, but in their paper, \cite{Anca2008} proves that none of these cases arise within the search space of the algorithm.  The search space for a sawtooth map which has comparable kneading data with the given map is the sawtooth map parameter isentrope
$$
L_{\mathcal{S}}(h^*) = \{(a,b) \in [\tfrac{1}{2}, 2]^2 : \log a + \log b = h^*\}.
$$
We check the kneading data comparability of any map within this space by checking the ordering of the d'itineraries $(q_\lambda,q_\mu)$ of the critical point $c_1 \in I_1$ of $q_\lambda$ and $c_2 \in I_2$ of $q_\mu$.  Firstly, note that the isentrope $L_{\mathcal{S}}(h^*)$ is actually a one dimensional curve, with a left boundary $(e^{h^* - \log 2}, 2)$ and right boundary $(2, e^{h^* - \log 2})$. Secondly, since $q_\lambda$ and $q_\mu$ are logistic maps, it is easy to show using differentiation that its only critical point is exactly $c_1 = c_2 = \frac{1}{2}$, and that in fact $I_1 = I_2$, but we distinguish them in order to make it clear whether the diorbits start $(x,f_1(x),f_2(f_1(x)),...)$ or $(x,f_2(x),f_1(f_2(x)),...)$.  We have the following lemma. (lemma 3.3 in \cite{Anca2008}).
\begin{lemma}
Let $(f_1,f_2)$ and $(f_1',f_2')$ be two distinct pairs of unimodal boundary anchored maps, and let $i(f,g)(x)$ denote the d'itinerary of $x$ under $(f,g)$.  If $i(f_1,f_2)(c_1) \leq i(f_1',f_2')(c_1)$ and $i(f_1,f_2)(c_2) \leq i(f_1',f_2')(c_2)$, then we have that $\bm{\mathcal{K}}(f_2 \circ f_1) \leq \bm{\mathcal{K}}(f_2' \circ f_1')$.
\end{lemma}

Radulescu's algorithm \cite{Anca2008} consists of two nested searches.  

The outer search is on the space of potential values of topological entropy that our given function, say, $f = q_{\mu} \circ q_{\lambda}$, could take.  We begin by setting a lower and upper bound for the entropy, then we progress by searching the isentrope $h^* = \frac{h_0 + h_1}{2}$ for a sawtooth map with comparable kneading data.  This midpoint then becomes one of the new bounds, and we repeat.  Now we know from earlier that the entropy of $f$ is contained in $[0,\log 4]$, however we must be careful, as it is not guaranteed that every sawtooth map isentrope we search is going to contain a sawtooth map with comparable kneading data.  However, from Corollay 6.1 in \cite{Anca2008} we know that
if we set the lower bound as $h_0 = 0$ and upper bound as $h_1 = \log(4.1)$, then our search will always eventually give us a map with comparable kneading data.

The inner search of the algorithm searches a specific isentrope for a map which has comparable kneading data with $f$.  Much like with the outer search, we check left boundary of the isentrope $L = (e^{h^* - \log 2}, 2)$, and the right boundary $R = (2, e^{h^* - \log 2})$.  If neither of these are comparable, then we check the midpoint $M = \frac{L+R}{2}$.  The sawtooth map with these parameters will either have kneading data comparable to that of $f$, or it will be a better upper or lower bound on the parameters of the map which does have comparable kneading data.  In particular, let $M = (M_1,M_2)$.  If 
\begin{align}
i(T_{M_1},T_{M_2})(c_1) &\leq i(q_\lambda,q_\mu)(c_1), \text{ and } \\
i(T_{M_1},T_{M_2})(c_2) &\geq i(q_\lambda,q_\mu)(c_2),
\end{align}
then $M$ is a better lower bound on the comparable parameter pair, and if
\begin{align}
i(T_{M_1},T_{M_2})(c_1) &\geq i(q_\lambda,q_\mu)(c_1), \text{ and } \\
i(T_{M_1},T_{M_2})(c_2) &\leq i(q_\lambda,q_\mu)(c_2),
\end{align}
then it is a better upper bound.  If it is a better lower bound, we reassign $L \leftarrow M$, and if it is a better upper bound, we reassign $R \leftarrow M$, and we repeat with a new $M \leftarrow \frac{L+R}{2}$.  When we find a map on the $h$-isentrope with comparable kneading data, say $S$ is the sawtooth map in question, then either:
\begin{equation}
\bm{\mathcal{K}}(S) \leq \bm{\mathcal{K}}(f),
\end{equation}
or
\begin{equation}
\bm{\mathcal{K}}(S) \geq \bm{\mathcal{K}}(f).
\end{equation}
In the first, we have found that $h$ is a better lower bound on the entropy of $f$, so we reassign $h_0 \leftarrow h$; in the second case, $h$ is a better upper bound on the entropy of $f$, so we reassign $h_1 \leftarrow h$.  We then repeat the above with a new $h \leftarrow \frac{h_0+h_1}{2}$ and continue in this way until $|h_1 - h_0| < \epsilon$ for some small $\epsilon > 0$.  In our implementation, we go by Radulescu's value of $\epsilon = 10^{-4}$.  In this implementation of the algorithm, we then take the average of the resulting lower and upper bound, and use this as our approximation of the entropy of $f$.  The final result of this implementation is the drawing of the isentropes in Figure 1.  For more details on the algorithm itself, see \cite{Anca2008}.

\begin{figure}[h!]
\centering
\includegraphics[width = 10cm]{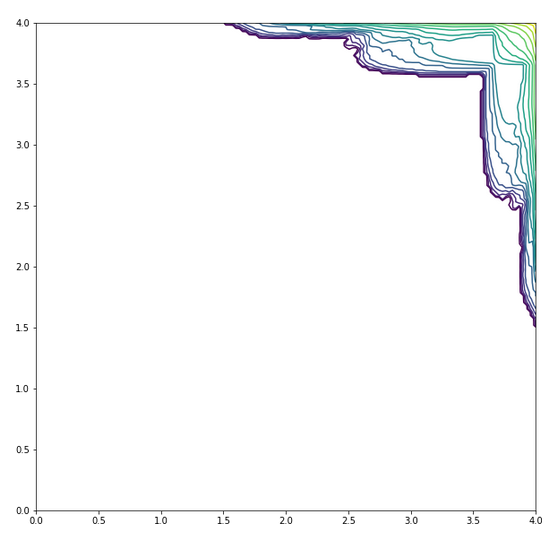}
\caption{\small{The isentropes of $P^Q$ in the parameter square $[0,4]^2$, where the point $(\lambda,\mu) \in [0,4]^2$ represents the map $f = q_\mu \circ q_\lambda$.  To compute entropy values in this space, we truncated the square $[0,4]^2$ into a $151 \times 151$ grid of points.}}
\end{figure}

As a final note on this algorithm, we see that the general framework for this algorithm is two nested bisection searches. The outer search is along the possible entropy values in the closed interval $[0,\log(m+1)]$ --  perhaps with a small perturbation of the upper bound so it becomes $[0, \log (m+1.1)]$ -- and an inner search along each isentrope in some model space, say tent maps, sawtooth maps, or even stunted sawtooth maps, for a map which has kneading data that is comparable with the original map.  The inner search tells us what ordering the entropy of the maps in the isentrope has.  The outer search provides the next candidate for the entropy of the map.  Together, they find progressively better bounds on the entropy of the map in question.  This framework could, with some extra work, provide specialised algorithms for many types of maps.  \cite{Block1992} is an example of this framework for cubic polynomials, using boundary anchored 2-modal sawtooth maps as the model mappings.  

One way that one could improve an algorithm like this is to refine the searching algorithms used in the nested searches.  Radulescu uses a pair of nested bisection searches, however, it might be faster to use interpolation searches.

\subsection{Amigo and Gimenez algorithm}

Now, we look at an algorithm by Amigo and Gimenez \cite{amigo2014} for computing general boundary anchored multimodal maps, based on an earlier algorithm by them and Rui Dil\~{a}o \cite{amigo2012}.   All the theoretical results in this section, unless stated otherwise, are from \cite{amigo2014}.

This algorithm seeks to calculate the lap number of $\ell( f^n )$ of  $f^n$  for each $n$, until the difference between consecutive steps is lower than some resolution.  I.e. given an $\varepsilon >0$, the algorithm stops when
$$
\Bigg| \frac{1}{N+1}\log \ell(f^{N+1}) - \frac{1}{N}\log \ell(f^N) \Bigg| < \varepsilon.
$$
For practical purposes, we also introduce another stopping criterion, that is given an integer $n_{max} \geq 2$, the algorithm stops if it has performed $N = n_{max}$ steps.  This prevents run times for the program becoming too long, and stops any potential infinite loops.  Of course, the natural question here is, given an $l$-modal boundary anchored map $f$, how do we calculate $\ell( f^n )$ for some positive integer $n$?  One answer of course would be to take a more general route, to calculate the solutions to 
$$
\frac{d f^n}{d x}(x) = 0,
$$
i.e. directly finding the critical points using numerical root-finding methods, and using this to deduce the lap number for $f^n$.  However, for large $n$ this can be slow, and many methods might risk missing critical points if the search is not careful enough.  Moreover, as we are calculating $\ell( f^n)$ for each $n$ up to a certain point, root-finding methods for finding $\ell( f^{N+1} )$ ignores all the computations we did to calculate $\ell( f^k )$ for $k < N$, thus making these methods inefficient in this context.  Instead, Amigo \textit{et al.} created an algorithm which uses the computations done in previous steps to aid the current step, by exploiting the property that $f^{N+1} = f \circ f^N$.  Specifically, it does this by firstly calculating only the critical points of $f$, and then studying the orbits of these critical points, in other words, studying the kneading data of $f$.
\par
The most important thing we will need from the orbits of the critical points is whether each element is a maximum or minimum, as this will give us information on the number of times $f^n$ turns.
 Fortunately, we have the following theorem.  For all of the results in this section, given an $l$-modal map $f$ of an interval $I = [a,b]$ with critical points $c_k$, we define the following partition of $I$ by $f$ as
$$
I_1 = [a,c_1), C_1 = \{ c_1 \}, I_2 = (c_1, c_2), ..., C_l = \{ c_l \}, I_{l+1} = (c_l, b].
$$

\begin{theorem}
\label{MinOrMax}
Let $f : I \rightarrow I$ be an $m$-modal map with positive shape, for some interval $I$, and let $n$ be a positive integer.  We have:
\begin{equation}
f^{n+1}(x) \hspace{3mm} \text{is a maximum if} \hspace{3mm} \left\{ \begin{aligned}
									                    &f^n(x) = c_{odd}, \\
										         &f^n(x) \in I_{even} \text{ and } f^n(x) \text{ is a minimum, or} \\
										         &f^n(x) \in I_{odd} \text{ and } f^n(x) \text{ is a maximum.}
                                                                                                                    \end{aligned} \right.
\end{equation}
Similarly,
\begin{equation}
f^{n+1}(x) \hspace{3mm} \text{is a minimum if} \hspace{3mm} \left\{ \begin{aligned}
									                   &f^n(x) = c_{even}, \\
										        &f^n(x) \in I_{odd} \text{ and } f^n(x) \text{ is a minimum, or} \\
										        &f^n(x) \in I_{even} \text{ and } f^n(x) \text{ is a maximum.}
                                                                                                                   \end{aligned} \right.
\end{equation}
\end{theorem}

Note that if the map has negative shape, then in (2.7), we replace ``$f^{n+1}(x)$ is a maximum if...'' with ``$f^{n+1}(x)$ is a minimum if...'', and the other way around in (2.8).  Using this, we will now define a expansion of the notion of kneading sequences.

\begin{definition}
We define the \textit{min-max sequences} of an $l$-modal map $f$ to be the vector $\bm{\omega} = (\omega_1, ..., \omega_l)$, where $\omega_i = (\omega_i^{(n)})_n$ is a sequence such that
$$
\omega_i^{(n)} = \left\{ \begin{aligned}
			         &m^{\mathcal{K}_i^{(n)}} \text{ if } f^n(c_i) \text{ is a minimum} \\
			         &M^{\mathcal{K}_i^{(n)}} \text{ if } f^n(c_i) \text{ is a maximum}
			         \end{aligned} \right. .
$$
\end{definition}

We see that the min-max sequence not only contain all the information that kneading sequences give, but also provide the extra information on whether $f^n(c_i)$ is a minimum or maximum.  Given the kneading sequences of $f$, we can use Theorem \ref{MinOrMax} to easily determine the min-max sequences of $f$.  Amigo then identifies two types of min max symbols with respect to a given critical line (that is, the vertical line $y = c_i$), called \textit{good symbols} and \textit{bad symbols}.  The good symbols and bad symbols will tell us how many critical points there are in $f^{n+1}$, the presence of good symbols indicating more, and the presence of bad symbols indicating fewer.  Essentially, the bad symbols with respect to the line $y = c_i$ are those min-max symbols that represent minima which are above the line $y = c_i$, and maxima which are below the line $y = c_i$.  Thus, bad symbols are the cases where critical points of $f^n$ don't 'produce' any more critical points in $f^{n+1}$, and so act as dead-ends.  This fortunately has an easy mathematical description - the set of bad symbols with respect to the line $y = c_i$ is

\begin{equation}
\mathcal{B}_i = \{ M^{I_1}, M^{C_1}, ..., M^{I_i}, M^{C_i}, m^{C_i}, m^{I_{i+1}}, ..., m^{C_l}, m^{I_{l+1}} \}.
\end{equation}

For $1 \leq i \leq l$ and $n \geq 0$, we let $s_i^{(n)} = |\{x \in I: f^n(x)- c_i = 0, f^k(x)-c_i \neq 0$ for $0 \leq k < n \}|$, i.e. the number of intersections between the curve $y = f^n(x)$ and the vertical line $y = c_i$, such that $y = c_i$ which don't intersect with $y = f^k(x)$ for any $0 \leq k < n$.  Note that $s_i^{(0)} = 1$ for all $i$.  Additionally, we set

\begin{equation}
\label{smalls}
s^{(n)} = \sum_{i = 1}^l s_i^{(n)},
\end{equation}

i.e. the total number of intersections with critical points.  We have the first part of the following theorem from \cite{amigo2012}.

\begin{theorem}
Let $\bm{\omega} = (\omega_1, ..., \omega_l)$ be the min-max sequences of an $l$-modal map $f$.  Then

\begin{equation}
\label{NewLap}
\ell( f^n ) = 1 + \sum_{k = 0}^{n - 1} s^{(k)}.
\end{equation}
\end{theorem}

The property of $s^{(k)}$ that the counting doesn't include points of intersection between $y = f^k(x)$ and any previous iteration of $f$ ($f^i(x)$ for $i < k$) essentially guarantees that when we perform the sum of thee $s^{(k)}$, we are not overcounting the number of intersections with a critical point as we iterate $f$.  Intuitively, this equation comes from the fact that the critical points are the maxima and minima of $f$, so if $f^k(x) = c_i$, then $f^{k+1}(x)$ is a maximum or minimum.  Hence, by counting $s^{(k)}$ up to $s^{(n-1)}$, we are counting the number of critical points of $f^n$, which means the minimum number of monotone pieces of $f^n$ must be one plus this number.  So now we have an equation for the lap number of $f^n$, we now wish to compute $s^{(k)}$, and by extension $s_i^{(k)}$.  In order to do this, we must count the $\omega_j^{(k)}$ which are bad symbols, $\mathcal{B}_i$.  Hence we let

$$
K_i^{(n)} = \{ (k,j): 1 \leq k \leq l, 1 \leq j \leq n, \omega_k^j \in \mathcal{B}_i \},
$$

which we can define in a recursive fashion by

\begin{equation}
\label{BadSet}
K_i^{(n)} \setminus K_i^{(n-1)} = \{ (k,n) : 1 \leq k \leq l, \omega_k^{(n)} \in \mathcal{B}_i \}.
\end{equation}

This formulation will be easier to implement in a program.  Using this set of indices, we define

\begin{equation}
\label{MissingPoints}
S_i^{(n)} = 2 \sum_{(k,j) \in K_i^{(n)}} s_k^{(n - j)},
\end{equation}

to be twice the number of times a minimum or maximum of $y = f^n(x)$ is above or below the line $y = c_i$, respectively.  Similarly to Equation \ref{smalls}, we set

\begin{equation}
\label{BigS}
S^{(n)} = \sum_{i = 1}^l S_i^{(n)}.
\end{equation}

\begin{figure}
\begin{subfigure}{.5\textwidth}
  \centering
  \includegraphics[width = 3in]{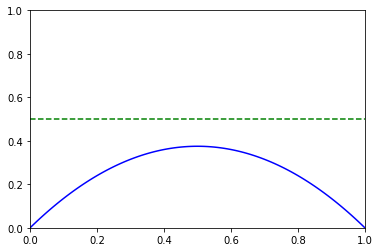}
\end{subfigure}%
\begin{subfigure}{.5\textwidth}
  \centering
  \includegraphics[width= 3in]{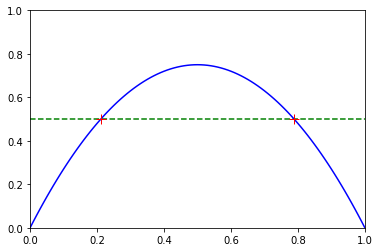}
\end{subfigure}
\caption{We see on these two figures that when the maximum of the function $f_a(x) = ax(1-x)$ is below the line $y = c = 0.5$, we don't get the two intersections we would get otherwise.  This is why we subtract twice the sum in equation \ref{MissingPoints}.}
\end{figure}

Bringing the above together, Theorem 2 from \cite{amigo2014} gives us a recursive definition for the $s_n^{(i)}$ values, which will give us the mathematical machinery we need to calculate the lap number of $f^n$.

\begin{theorem}
Let $f$ be an $l$-modal boundary anchored map.  Then
\begin{equation}
\label{Iteration}
s_i^{(n)} = 1 + \sum_{k = 0}^{n - 1} s^{(k)} - S_i^{(n)}.
\end{equation}
\end{theorem}

This finally gives us a nice formulation for the lap number of $f^n$:
\begin{equation}
\label{Laps}
\ell( f^n ) = \frac{s^{(n)} + S^{(n)}}{l}.
\end{equation}

With this final result, we move on to describing this algorithm for computing topological entropy.

The first step to perform in this algorithm is to start off our calculations.  After calculating the min-max sequences of $f$, then for $1 \leq i \leq l$, we set:
\begin{align*}
&s_i^{(0)} = 1, \text{ and } \\
&K_i^{(1)} = \{(k,1) : 1 \leq k \leq l, \omega_k^{(1)} \in \mathcal{B}_i\}.
\end{align*}

Using these, we initialise the rest of the $n=1$ variables as follows:

\begin{align*}
&S_i^{(1)} = 2 \sum_{(k,1) \in K_i^{(1)}} s_0^{(k)}, \\
&S^{(1)} = \sum_{i = 1}^l S_i^{(1)}, \\
&s_i^{(1)} = 1 + s^{(0)} - S_i^{(1)} = 1 + l - S_i^{(1)}, \text{ and } \\
&s^{(1)} = l(1 + s^{(0)}) - S^{(1)} = l^2 + l - S^{(1)}.
\end{align*}

Now, we are ready to begin the general loop of the algorithm.  For $n \geq 2$, and for all $1 \leq i \leq l$, we recursively calculate the following, in order of presentation:

\begin{align*}
&K_i^{(n)} \text{ from } K_i^{(n-1)} \text{ using equation \ref{BadSet}}, \\
&S_i^{(n)} \text{ from } s_i^{(0)}, ..., s_i^{(n-1)}, K_i^{(n)} \text{ using equation \ref{MissingPoints}}, \\
&S^{(n)} \text{ from } S_i^{(n)} \text{ using equation \ref{BigS}}, \\
&s_i^{(n)} \text{ from } s_i^{(0)}, ..., s_i^{(n-1)}, S_i^{(n)} \text{ using equation \ref{Iteration}}, \\
&s^{(n)} \text{ from } s_i^{(n)} \text{ using equation \ref{smalls}}, \\
&\ell( f^n ) \text{ from } s^{(n)}, S^{(n)} \text{ using equation \ref{Laps}}.
\end{align*}

We terminate this loop if one of two cases occur.  In the first case, we have that

\begin{equation}
\Bigg| \frac{1}{n} \log \ell( f^n ) - \frac{1}{n - 1} \log \ell( f^{n-1} ) \Bigg| < \varepsilon,
\end{equation}

for some pre-chosen resolution $\varepsilon$, in which case we stop and output

\begin{equation}
h(f) = \frac{1}{n} \log \ell( f^n ) = \frac{1}{n} \log \frac{s^{(n)} + S^{(n)}}{l}.
\end{equation}

In the second case, we have that $n = n_{max}$, in which case the algorithm has failed.
\par
In Appendix B, the code for our implementation can be seen, in which we have chosen the resolution to be $\varepsilon = 10^{-4}$, and $n_{max} = 2000$.

\subsection{Comparing the two algorithms}

Radulescu's algorithm \cite{Anca2008} works by finding a model for the function in question which has similar enough entropy, and thus requires the usage of search algorithms, with the aid of partial orders, to find the entropy of the map.  On the other hand, Amigo and Giminez's algorithm \cite{amigo2014} studies the iterations of the given function, specifically with regards to orbits of the critical points, i.e. if they are maxima or minima, as well as the kneading data.  Using this, the algorithm performs arithmetic to count the number of critical points there are in each iteration, compute the lap number from this, and thus compute the topological entropy.   

The benefit of using an algorithm like Radulescu's is that it exploits the easy formulation of the topological entropy of sawtooth maps, which lends to its efficiency.  
The main problem, however, with this algorithmic structure is that it relies heavily on the accuracy of the orbits of the critical points, since the algorithms directly compare kneading data, so an implementation of an algorithm like this may have to resort to variable precision arithmetic to ensure the kneading data are correct - if the kneading data are incorrect, there is no way to guarantee that there will exist a sawtooth map in a particular isentrope which has kneading data that is comparable to the incorrect kneading data for the map $f$.  For the case of Radulescu's algorithm, this isn't as much of a problem, since the critical point is always $0.5$, a rational number with a short decimal expansion, and the shape of the logistic function is simpler than that of a higher degree polynomial.  However, for higher degree polynomials, numerical methods may need to be used to compute the critical points.  Hence, more computation power and time is required to get the accuracy needed for computing entropy this way.

The benefit to using Amigo's algorithm is, without needing to change any of the algorithm, it works for any interval map, so long as it is boundary anchored and possesses a finite number of critical points.  For instance, we could find the entropy of $f_\lambda : [0,\pi] \rightarrow [0,\pi]; x \mapsto \lambda \sin(x)$ for some $\lambda$ with $0 < \lambda \leq \pi$.  

One problem with Amigo's algorithm is that it is often very slow to converge for higher precision.  In their \cite{amigo2014} paper, they explain that the resolution $\varepsilon$ for the stopping criterion doesn't bound the error of the approximation with respect to the true value of entropy, instead it bounds the error between consecutive iterations of the algorithm.  This implies, and shows in numerical simulations later in that paper, that in order to achieve high precision for the value of entropy, one may need $\varepsilon$ far smaller than the desired accuracy.  This comes with the price that to achieve this precision, one needs to perform very many iterations of the algorithm.  For instance, we have found that computing the entropy of some cubic polynomials of negative shape to a resolution of $\varepsilon = 10^{-4}$ needs more than $1500$ iterations.

\section{Topological entropy of boundary anchored cubic polynomials}

In this section we implement the algorithm in \cite{amigo2014}.
The first step to perform in this algorithm is to start off our calculations.  After calculating the min-max sequences of $f$, then for $1 \leq i \leq l$, we set:
\begin{align*}
&s_i^{(0)} = 1, \text{ and } \\
&K_i^{(1)} = \{(k,1) : 1 \leq k \leq l, \omega_k^{(1)} \in \mathcal{B}_i\}.
\end{align*}

Using these, we initialise the rest of the $n=1$ variables as follows:

\begin{align*}
&S_i^{(1)} = 2 \sum_{(k,1) \in K_i^{(1)}} s_0^{(k)}, \\
&S^{(1)} = \sum_{i = 1}^l S_i^{(1)}, \\
&s_i^{(1)} = 1 + s^{(0)} - S_i^{(1)} = 1 + l - S_i^{(1)}, \text{ and } \\
&s^{(1)} = l(1 + s^{(0)}) - S^{(1)} = l^2 + l - S^{(1)}.
\end{align*}

Now, we are ready to begin the general loop of the algorithm.  For $n \geq 2$, and for all $1 \leq i \leq l$, we recursively calculate the following, in order of presentation:

\begin{align*}
&K_i^{(n)} \text{ from } K_i^{(n-1)} \text{ using equation \ref{BadSet}}, \\
&S_i^{(n)} \text{ from } s_i^{(0)}, ..., s_i^{(n-1)}, K_i^{(n)} \text{ using equation \ref{MissingPoints}}, \\
&S^{(n)} \text{ from } S_i^{(n)} \text{ using equation \ref{BigS}}, \\
&s_i^{(n)} \text{ from } s_i^{(0)}, ..., s_i^{(n-1)}, S_i^{(n)} \text{ using equation \ref{Iteration}}, \\
&s^{(n)} \text{ from } s_i^{(n)} \text{ using equation \ref{smalls}}, \\
&\ell( f^n ) \text{ from } s^{(n)}, S^{(n)} \text{ using equation \ref{Laps}}.
\end{align*}

We terminate this loop if one of two cases occur.  In the first case, we have that

\begin{equation}
\Bigg| \frac{1}{n} \log \ell( f^n ) - \frac{1}{n - 1} \log \ell( f^{n-1} ) \Bigg| < \varepsilon,
\end{equation}

for some pre-chosen resolution $\varepsilon$, in which case we stop and output

\begin{equation}
h(f) = \frac{1}{n} \log \ell( f^n ) = \frac{1}{n} \log \frac{s^{(n)} + S^{(n)}}{l}.
\end{equation}

In the second case, we have that $n = n_{max}$, in which case the algorithm has failed.

\begin{figure}
\begin{subfigure}{.5\textwidth}
  \centering
  \includegraphics[width = 3in]{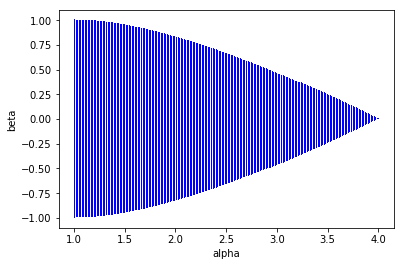}
\end{subfigure}%
\begin{subfigure}{.5\textwidth}
  \centering
  \includegraphics[width= 3in]{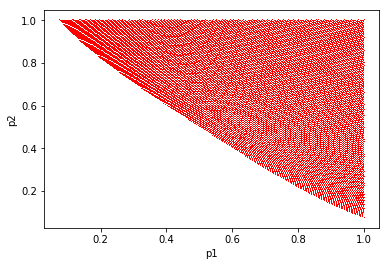}
\end{subfigure}
\caption{Here we see a lattice of grid points which show the $(\alpha,\beta)$ parameter space and the $(p_1,p_2)$ parameter space, respectively.}
\end{figure}

Now, in order to compute the topological entropy of the boundary anchored cubic polynomials, we must find a suitable parameterisation of them which preserves the property of them being boundary anchored.  Fortunately, we can derive such a parameterisation into two parameters as follows.  Assuming $I = [-1,1]$ and $\sigma_1 \in \{-1,1\}$, let

$$
f(x) = a x^3 + b x^2 + c x + d,
$$

then we have the following conditions

\begin{align*}
&f(-1) = -a+b-c+d = \sigma_1, \\
&f(1) = a+b+c+d = -\sigma_1
\end{align*}

By adding and subtracting the two equations, we get a pair of relations

\begin{align*}
b+d = 0, \\
a+c = -\sigma_1,
\end{align*}

and so by setting $a = \alpha$ and $b = \beta$, we get

\begin{align*}
d = -\beta, \\
c = -\sigma_1 - \alpha.
\end{align*}

This gives the formulation of $f$ as

\begin{equation}
f(x) = \alpha x^3 + \beta x^2 + (-\sigma_1 - \alpha) x - \beta.
\end{equation}

From \cite{bruin2013}, we have the ranges of $\alpha$ and $\beta$ for cubic polynomials of positive shape, for which $f(I) \subseteq I$, are $0 \leq \alpha \leq 4$ and $|\beta| \leq 2 \sqrt{\alpha} - \alpha$.  Perfoming some tests reveals that the critical points are only real and distinct iff $\alpha > 1$, so it will be sufficient to consider the range $1 < \alpha \leq 4$.  In a similar manner, it is easy to see from numerical simulations that a similar ranges of $\alpha$ and $\beta$ for cubic polynomials of negative shape, for which the critical values are real and distinct, and $f(I) \subseteq I$, are $-4 \leq \alpha < -1$ and $|\beta| < 2 \sqrt{|\alpha|} - |\alpha|$.  Now that we have the parameterisation and the ranges for the parameters, we will compute the critical points as follows:

\begin{align*}
c_{1,2} = -\frac{\beta}{3\alpha} \pm \frac{1}{3\alpha}\sqrt{\beta^2 - 3\alpha(1 - \alpha)},
\end{align*}

which is trivially derived from solving $f'(x) = 0$. From this, we can compute the critical values, and use the resultant parameter pair $(p_1,p_2)$ as the coordinates to which we assign the topological entropy of $f$.  This way, we don't need to find a direct parameterisation of $f$ in terms of the critical values.  Using this, we produce contour plots for the topological entropy of the positive and negative shape, boundary anchored cubic polynomials.
\par

\begin{figure}[h!]
\centering
\includegraphics[width = 10cm]{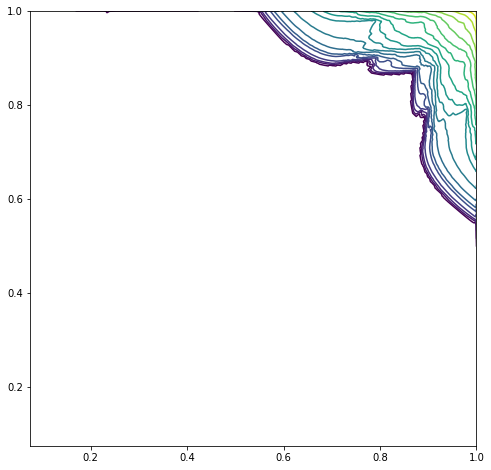}
\caption{This is the contour plot of the topological entropy of the positive shape, boundary anchored cubic polynomials.}
\label{PositiveShape}
\end{figure}

\begin{figure}[h!]
\centering
\includegraphics[width= 10cm]{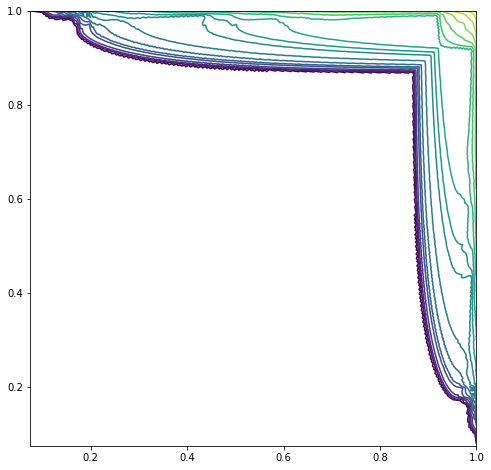}
\caption{This is the contour plot of the topological entropy of the negative shape boundary anchored cubic polynomials.}
\label{NegativeShape}
\end{figure}

We see from the plots on figures \ref{PositiveShape} and \ref{NegativeShape} that the dynamics of positive shape cubic polynomials differs significantly from that of negative shape cubics.  Looking at the shape of the isentropes near the boundary between 0 and positive topological entropy in figure \ref{PositiveShape}, it seems that some isentropes may not be locally connected around the central area, where it seems the shape becomes rougher.  If this roughness continues indefinitely as one zooms in, or if other strange behaviour occurs here, then it is possible that isentropes which pass through that area may not be locally connected.
\par
In figure \ref{NegativeShape}, the local connectedness is questionable around the tips of the two branches, where the isentropes appear to bunch together and become craggly in shape.  Again, if this craggliness continues indefinitely as one zooms in, or some other strange behaviour occurs here, then it is possible that isentropes that pass through this area may not be locally connected.
\par
Of course, the contours that appear on these plots only show some of the isentropes.  To gain an even better insight into the shape of the surface of topological entropy, more contours can be added to the function which plots the contours, and for more powerful computers, one can make the resolution go down to $\varepsilon = 10^{-6}$ or even smaller to get more accuracy, though one must remember to also increase $n_{max}$ accordingly.  Moreover, the number of points used in the grid lattice of the parameter space can be increased by changing the 'Precision' argument in the draw functions provided, thus allowing more detail to be shown.  Additionally, with some extra programming, one can alter the draw functions to plot the isentropes on a particular subset of the parameter triangle $P^2$, to effectively zoom in on the isentropes to view the details clearer.

\textbf{Acknowledgements.} The authors would like to thank Sebastian van Strien and Trevor Clark for valuable suggestions and discussions.

\end{document}